\documentclass[12pt]{article}%
\usepackage{amsfonts}
\usepackage{amsmath,amssymb}
\usepackage[applemac]{inputenc}
\usepackage{amsmath,amssymb,fullpage}
\usepackage{color}
\usepackage{showlabels}
\usepackage{amsmath}
\usepackage{amssymb}
\usepackage{graphicx}%
\providecommand{\U}[1]{\protect\rule{.1in}{.1in}}

\newcommand{\dproof}{\noindent {Proof.} \quad}
\newcommand{\fproof}{\hfill $\square$ \bigskip}
\newtheorem{definition}{Definition}[section]
\newtheorem{example}{Example}[section]
\newtheorem{theorem}[definition]{Theorem}
\newtheorem{problem}[definition]{Problem}
\newtheorem{remark}[definition]{ \it Remark}

\newtheorem{proposition}[definition]{Proposition}
\newtheorem{lemma}[definition]{Lemma}

\numberwithin{equation}{section}

\def\EE{{\mathbb{ E}}}

\def\1B{\text{1\!\!I}}

\begin{document}

\date{9 October 2018}
\title{Singular control and optimal stopping of memory mean-field processes}
\author{Nacira Agram$^{1,2}$, Achref Bachouch$^{1}$, Bernt \O ksendal$^{1}$ and Frank
Proske$^{1}$ }
\maketitle


\footnotetext[1]{Department of Mathematics, University of Oslo, P.O. Box 1053
Blindern, N--0316 Oslo, Norway.  Emails: \texttt{naciraa@math.uio.no, achrefb@math.uio.no,
oksendal@math.uio.no, proske@math.uio.no}.\newline This research was carried
out with support of the Norwegian Research Council, within the research
project Challenges in Stochastic Control, Information and Applications
(STOCONINF), project number 250768/F20.}
\footnotetext[2]{ University Mohamed Khider, Biskra, Algeria.
}

\begin{abstract}
The purpose of this paper is to study the following topics and the relation between them:\\
(i) Optimal singular control of mean-field stochastic differential equations with memory,\\
(ii) reflected advanced mean-field backward stochastic differential equations, and\\
(iii) optimal stopping of mean-field stochastic differential equations.\\

More specifically, we do the following:

\begin{itemize}
\item We prove the existence and uniqueness of the solutions of some reflected advanced memory backward stochastic differential equations (AMBSDEs),

\item we give sufficient and necessary conditions for an optimal singular
control of a memory mean-field stochastic differential equation (MMSDE) with partial information, and

\item we deduce a relation between the optimal singular control of a MMSDE, and the optimal stopping of such processes.

\end{itemize}
\end{abstract}

\section{Introduction}

\noindent Let $(%
\Omega
,\mathcal{F},\mathbb{P})$ be a given probability space with filtration
$\mathbb{F}=(\mathcal{F}_{t})_{t\geq0}$ generated by a 1-dimensional Brownian
motion $B=B(t,\omega); (t,\omega) \in[0,T] \times\Omega.$ Let $\mathbb{G}%
=\{\mathcal{G}_{t}\}_{t\geq0}$ be a given subfiltration of $\mathbb{F}%
=(\mathcal{F}_{t})_{t\geq0}$ , in the sense that $\mathcal{G}_{t}%
\subset\mathcal{F}_{t}$ for all $t.$

\noindent The purpose of this paper is to study the following concepts and problems, and
the relation between them. For simplicity of notation we deal only with the
$1$-dimensional case.\newline

\begin{itemize}
\item \emph{Topic 1: Optimal singular control of memory mean-field stochastic
differential equations:}\newline

Consider the following \emph{mean-field memory singular controlled system},
with a state process $X(t)=X^{\xi}(t)$ and a singular control process
$\xi(t),$ of the form
\begin{equation}
\left\{
\begin{array}
[c]{l}%
dX(t)=b(t,X(t),X_{t},M(t),M_{t},\xi(t),\omega)dt+\sigma(t,X(t),X_{t}%
,M(t),M_{t},\xi(t),\omega)dB(t)\\
\quad\quad\quad+\lambda(t,\omega)d\xi(t);\quad t\in\lbrack0,T],\\
X(t)=\alpha(t);\quad t\in\lbrack-\delta,0],
\end{array}
\right.  \label{eq6.1a}%
\end{equation}
where
\begin{align*}
&  X_{t}=\{X(t-s)\}_{0\leq s\leq\delta},\quad\text{(the \emph{memory segment
of}}X(t)),\\
&  M(t)=\mathcal{L}(X(t))\quad\text{(the \emph{law of}}X(t)),\\
&  M_{t}=\{M(t-s)\}_{0\leq s\leq\delta},\quad\text{(the \emph{memory segment
of}}M(t)).
\end{align*}
We assume that our control process $\xi(t)$ is ${\mathbb{R}}$-valued
right-continuous $\mathbb{G}$-adapted process, and $t\mapsto\xi(t)$ is
increasing (non-decreasing) with $\xi(0^{-})=0$, and such that the
corresponding state equation has a unique solution $X$ with $\omega\mapsto
X(t,\omega)\in L^{2}(\mathbb{P})$ for all $t$. The set of such processes $\xi$
is denoted by $\Xi$. \newline


The \emph{performance functional} is assumed to be of the form
\begin{align*}
J(\xi) &  ={\mathbb{E}}[%
{\textstyle\int_{0}^{T}}
f(t,X(t),X_{t},M(t),M_{t},\xi(t),\omega)dt+g(X(T),M(T),\omega)\\
&  \qquad+%
{\textstyle\int_{0}^{T}}
h(t,X(t),\omega)d\xi(t)];\quad\xi\in\Xi\,.
\end{align*}
For simplicity we will in the following suppress the $\omega$ in the notation.
\newline We may interpret these terms as follows:\newline The state $X(t)$ may
be regarded as the value at time $t$ of, e.g. a fish population. The control
process $\xi(t)$ models the amount harvested up to time $t$, the coefficient
$\lambda(t)$ is the unit price of the amount harvested, $f$ is a profit rate,
$g$ is a bequest or salvage value function, and $h$ is a cost rate for the use
of the singular control $\xi$. The $\sigma$-algebra $\mathcal{G}_{t}$
represents the amount of information available to the controller at time $t.$
The problem we consider, is the following:

\begin{problem}
Find an optimal control $\hat{\xi}\in\Xi$ such that
\begin{equation}
J(\hat{\xi})=\sup_{\xi\in\Xi}J(\xi)\,.\label{eq6.4}%
\end{equation}

\end{problem}

This problem turns out to be closely related to the following topic:
\end{itemize}

\begin{itemize}
\item \emph{Topic 2: Reflected mean-field backward stochastic differential
equations}
\end{itemize}

We study reflected AMBSDEs where at any time $t$ the driver $F$
may depend on future information of the solution processes.
More precisely, for a given driver $F$, a given threshold process $S(t)$ and a
given terminal value $R$ we consider the following type of reflected AMBSDEs in the
unknown processes $Y,Z,K$:

\begin{equation}
\left\{
\begin{array}
[c]{l}%
(i)Y(t)=R+%
{\textstyle\int_{t}^{T}}
F(s,Y(s),Z(s),\mathbb{E}[Y^{s}|\mathcal{F}_{s}],\mathbb{E}[Z^{s}%
|\mathcal{F}_{s}],\mathcal{L}(Y^{s},Z^{s}))ds\\
\quad\quad\quad\quad\quad\quad+K(T)-K(t)-%
{\textstyle\int_{t}^{T}}
Z(s)dB(s); \quad0\leq t\leq T,\\
(ii)Y(t)\geq S(t); \quad0\leq t\leq T,\\
(iii){%
{\textstyle\int_{0}^{T}}
}(Y(t)-S(t))dK^{c}(t)=0\text{ a.s. and }\triangle K^{d}(t)=-\triangle
Y(t)\mathbf{1}_{\{Y(t^{-})=S(t^{-})\}}\text{ a.s.},\\
(iv)Y(t)=R; \quad t\geq T,\\
(v)Z(t)=0; \quad t>T.
\end{array}
\right.
\end{equation}
\vskip 0.3cm
\noindent Here $\mathcal{L}(Y^{s},Z^{s})$ is the joint law of paths $(Y^{s},Z^{s})$, and
for a given positive constant $\delta$ we have put
\[
\quad Y^{t}:=\{Y(t+s)\}_{s\in\lbrack0,\delta]}\text{ and }Z^{t}:=\{Z(t+s)\}_{s\in\lbrack0,\delta]}\text{ (the (time)-advanced
segment). }%
\]
\vskip0.2cm

This problem is connected to the following:
\begin{itemize}
\item \emph{Topic 3: Optimal stopping and its relation to the problems above.}
\end{itemize}
For $t \in [0,T]$ let $\mathcal{T}_{[t,T]}$ denote the set of all $\mathbb{F}$-stopping
times $\tau$ with values in $[t,T].$\\
Suppose $\left(  Y,Z,K\right)  $ is a solution of the reflected AMBSDE in Topic 2 above.
\begin{description}
\item[(i)] Then, for $t\in\left[  0,T\right]$, the process $Y(t)$ is the solution of the optimal stopping problem%
\begin{align}
Y(t)=\underset{\tau\in\mathcal{T}_{[t,T]}}{ess\sup}\ \Big\{ \mathbb{E} [%
{\textstyle\int_{t}^{\tau}}
&F(s,Y(s), Z(s), \mathbb{E}[Y^{s}|\mathcal{F}_s],\mathbb{E}[Z^{s}|\mathcal{F}_s],\mathcal{L}(Y^{s},Z^{s}))ds\nonumber\\
&+S(\tau)\mathbf{1}_{\tau<T}%
+R\mathbf{1}_{\tau=T}|\mathcal{F}_{t}]\Big\}.
\end{align}

\item[(ii)] Moreover, for $t \in [0,T]$  the solution process $K(t)$ is given by%
\begin{align}
&K(T)-K(T-t)\nonumber\\
&=\underset{s\leq t}{\max}\Big\{R+%
\int_{T-s}^{T}
F(r,Y(r), Z(r), \mathbb{E}[Y^{r}|\mathcal{F}_r],\mathbb{E}[Z^{r}|\mathcal{F}_r],\mathcal{L}(Y^{r},Z^{r}))dr\nonumber\\
&-\int _{T-s}^{T} Z(r)dB(r)-S(T-s)\Big\}^{-},
\end{align}
where $x^{-}=\max(-x,0),$
 and an optimal stopping time $\hat{\tau}_{t}$ is given by%
\begin{align*}
\hat{\tau}_{t}:&=\inf\{s\in\lbrack t,T],Y(s)\leq S(s)\}\wedge T\\
&  =\inf\{s\in\lbrack t,T],K(s) > K(t)\}\wedge T.
\end{align*}
\item[(iii)] In particular, if we choose $t=0$, we get that
\begin{align*}
\hat{\tau}_{0}:&=\inf\{s\in\lbrack 0,T],Y(s)\leq S(s)\}\wedge T\\
&  =\inf\{s\in\lbrack 0,T],K(s) > 0\}\wedge T
\end{align*}
solves the optimal stopping problem
\begin{align}
Y(0)&=\sup_{\tau\in\mathcal{T}_{[0,T]}}\mathbb{E}[{\textstyle\int_{0}^{\tau}}
F(s,Y(s), Z(s), \mathbb{E}[Y^{s}|\mathcal{F}_s],\mathbb{E}[Z^{s}|\mathcal{F}_s],\mathcal{L}(Y^{s},Z^{s}))ds\nonumber\\
&+S(\tau)\mathbf{1}_{\tau<T}%
+R\mathbf{1}_{\tau=T}],t\in\left[  0,T\right]  .
\end{align}

\end{description}

\noindent More specifically, the content of the paper is the
following:\newline

\noindent In Section 2, we define the spaces of measures and spaces of path
segments with their associated norms, and we give the necessary background
results for our methods.\newline

\noindent In Section 3, we prove existence and uniqueness of the solution for a
class of reflected advanced mean-field backward stochastic differential
equations.
\newline

\noindent In Section 4, we recall a fundamental connection between a
class of reflected AMBSDEs and optimal stopping under partial information.
equations.
\newline

\noindent Then in Section 5, we study the problem of optimal singular control of
memory mean-field stochastic differential equations. We give sufficient and necessary conditions for optimality in terms of variational inequalities.
\newline

\noindent Finally, in Section 6, we deduce a relation between the following quantities:\\
(i) The solution of a singular control problem for a mean-field SDE with memory. \\
(ii) The solution of a coupled system of forward memory \& backward advanced mean-field SDEs.\\
(iii) The solution of an optimal stopping problem involving these quantities.
\newline

\section{A Hilbert space of random measures}

\noindent In this section, we proceed as in Agram and \O ksendal \cite{AO3}, \cite{AO2} and construct a Hilbert space $\mathcal{M}$ of random measures on
$\mathbb{R}$. It is simpler to work with than the Wasserstein metric space
that has been used by many authors previously. See e.g. Carmona \textit{et al}
\cite{carmona1}, \cite{carmona}, Buckdahn \textit{et al} \cite{BLP} and the
references therein.

\noindent Following Agram and \O ksendal  \cite{AO3}, \cite{AO2}, we now
introduce the following Hilbert spaces:

\begin{definition}
\-

\begin{itemize}
\item Let $n$ be a given natural number. Then we define $\tilde{\mathcal{M}%
}=\tilde{\mathcal{M}}^{n}$ to be the pre-Hilbert space of random measures
$\mu$ on $\mathbb{R}^{n}$ equipped with the norm
\[%
\begin{array}
[c]{lll}%
\left\Vert \mu\right\Vert _{\tilde{\mathcal{M}}^{n}}^{2} & := & \mathbb{E[}%
{\textstyle\int_{\mathbb{R}^n}}
|\hat{\mu}(y)|^{2}(1+|y|)^{-2}dy]\text{,}%
\end{array}
\]
with $y=(y_{1},y_{2}, ... ,y_{n})\in\mathbb{R}^{n}$, and $\hat{\mu}$ is the
Fourier transform of the measure $\mu$, i.e.%
\[%
\begin{array}
[c]{lll}%
\hat{\mu}(y) & := & {%
{\textstyle\int_{\mathbb{R}^n}}
}e^{-ixy}d\mu(x);\quad y\in\mathbb{R}^{n},
\end{array}
\]
where $xy =x \cdot y = x_{1} y_{1} + x_{2} y_{2} + ... + x_{n} y_{n}$ is the
scalar product in $\mathbb{R}^{n}$.

\item $\tilde{\mathcal{M}}_{\delta}$ is the pre-Hilbert space of all path
segments $\overline{\mu}=\{\mu(s)\}_{s\in[0,\delta]}$ of processes $\mu
(\cdot)$ with $\mu(s)\in\tilde{\mathcal{M}}$ for each $s\in[0,\delta]$,
equipped with the norm
\begin{equation}
\left\Vert \overline{\mu}\right\Vert ^{2}_{\tilde{\mathcal{M}}_{\delta}}:={%
{\textstyle\int_0^{\delta}}
}\left\Vert \mu(s)\right\Vert ^{2}_{\tilde{\mathcal{M}}}ds.
\end{equation}

\item We let $\mathcal{M}$ and $\mathcal{M}_{\delta}$ denote the completion of
$\tilde{\mathcal{M}}$ and $\tilde{\mathcal{M}}_{\delta}$ and we let
$\mathcal{M}_{0}$ and $\mathcal{M}_{0,\delta}$ denote the set of deterministic
elements of $\mathcal{M}$ and $\mathcal{M}_{0,\delta}$, respectively.

\end{itemize}
\end{definition}

\noindent There are several advantages with working with this Hilbert space
$\mathcal{M}$, compared to the Wasserstein metric space:

\begin{itemize}
\item A Hilbert space has a useful stronger structure than a metric space.
\item Our space $\mathcal{M}$ is easier to work with.
\item The Wasserstein metric space $\mathcal{P}_{2}$ deals only with
probability measures with finite second moment, while our Hilbert space deals
with any (possibly random) measure $\mu\in\mathcal{M}$.
\end{itemize}

\noindent Let us give some examples for $n=1$:

\begin{example}
[Measures]\-

\begin{enumerate}
\item Suppose that $\mu=\delta_{x_{0}}$, the unit point mass at $x_{0}%
\in\mathbb{R}$. Then $\delta_{x_{0}}\in\mathcal{M}_{0}$ and
\[
{%
{\textstyle\int_{\mathbb{R}}}
}e^{ixy}d\mu(x)=e^{ix_{0}y},
\]

and hence
\[%
\begin{array}
[c]{lll}%
\left\Vert \mu\right\Vert _{\mathcal{M}_{0}}^{2} & =%
{\textstyle\int_{\mathbb{R}}}
|e^{ix_{0}y}|^{2}(1+|y|)^{-2}dy & <\infty\text{.}%
\end{array}
\]

\item Suppose $d\mu(x)=f(x)dx$, where $f\in L^{1}(\mathbb{R})$. Then $\mu
\in\mathcal{M}_{0}$ and by Riemann-Lebesque lemma, $\hat{\mu}(y)\in
C_{0}(\mathbb{R})$, i.e. $\hat{\mu}$ is continuous and $\hat{\mu
}(y)\rightarrow0$ when $|y|\rightarrow\infty$. In particular, $|\hat{\mu}|$ is
bounded on $\mathbb{R}$ and hence%
\[%
\begin{array}
[c]{lll}%
\left\Vert \mu\right\Vert _{\mathcal{M}_{0}}^{2} & =%
{\textstyle\int_{\mathbb{R}}}
|\hat{\mu}(y)|^{2}(1+|y[)^{-2}dy & <\infty\text{.}%
\end{array}
\]

\item Suppose that $\mu$ is any finite positive measure on $\mathbb{R}$. Then
$\mu\in\mathcal{M}_{0}$ and
\[%
\begin{array}
[c]{lll}%
|\hat{\mu}(y)| & \leq%
{\textstyle\int_{\mathbb{R}}}
d\mu(y)=\mu(\mathbb{R}) & <\infty\text{ for all }y\text{,}%
\end{array}
\]
and hence%
\[%
\begin{array}
[c]{lll}%
\left\Vert \mu\right\Vert _{\mathcal{M}_{0}}^{2} & =%
{\textstyle\int_{\mathbb{R}}}
|\hat{\mu}(y)|^{2}(1+|y|)^{-2}dy & <\infty\text{.}%
\end{array}
\]

\item Next, suppose $x_{0}=x_{0}(\omega)$ is random. Then $\delta
_{x_{0}(\omega)}$ is a random measure in $\mathcal{M}$. Similarly, if
$f(x)=f(x,\omega)$ is random, then $d\mu(x,\omega)=f(x,\omega)dx$ is a random
measure in $\mathcal{M}$.\newline
\end{enumerate}
\end{example}
\begin{definition}
[Law process]From now on we use the notation
\[
M_{t}:=M(t):=\mathcal{L}(X(t));\quad0\leq t\leq T,
\]
for the law process $\mathcal{L}(X(t))$ of $X(t)$ with respect to the probability $\mathbb{P}$.
\end{definition}
We recall the following results from Agram \& \O ksendal \cite{AO3}:

\begin{lemma}
\label{m'}The map $t\mapsto M(t):[0,T]\rightarrow\mathcal{M}_{0}$ is
absolutely continuous, and the derivative%

\[
M^{\prime}(t):=\frac{d}{dt}M(t)
\]
exists for all $t$.
\end{lemma}

\begin{lemma}
If $X(t)$ is an Itô-Lévy process as in \eqref{eq6.1a}, then the derivative
$M^{\prime}(s):=\frac{d}{ds}M(s)$ exists in $\mathcal{M}_{0}$ for a.a. $s$,
and we have
\[
M(t)=M(0)+%
{\textstyle\int_{0}^{t}}
M^{\prime}(s)ds;\quad t\geq0.
\]

\end{lemma}

The following result, based on Agram \& \O ksendal \cite{AO2},  is essential for our approach:

\begin{lemma}
\label{Lemma 4} \-

\begin{description}
\item[(i)] Let $X^{(1)}$ and $X^{(2)}$ be two $2$-dimensional random variables
in $L^{2}(\mathbb{P})$. Then there exist a constant $C_{0}$ not depending on
$X^{(1)}$ and $X^{(2)}$, such that
\[%
\begin{array}
[c]{lll}%
\left\Vert \mathcal{L}(X^{(1)})-\mathcal{L}(X^{(2)})\right\Vert _{\mathcal{M}%
_{0}^{2}}^{2} & \leq & C_{0}\ \mathbb{E}[(X^{(1)}-X^{(2)})^{2}]\text{.}%
\end{array}
\]

\item[(ii)] Let $\{X^{(1)}(t)\}_{t\in [0,T]},$ $\{X^{(2)}(t)\}_{t\in [0.T]}$ be two
paths, such that
\[
\mathbb{E}[%
{\textstyle\int_{0}^{T}}
X^{(i)2}(s)ds]<\infty \text{ for }i=1,2\text{.}%
\]
Then, for all $t$,
\[%
\begin{array}
[c]{lll}%
||\mathcal{L}(X_{t}^{(1)})-\mathcal{L}(X_{t}^{(2)})||_{\mathcal{M}_{0,\delta
}^{2}}^{2} & \leq & C_{0}\ \mathbb{E}[%
{\textstyle\int_{-\delta}^{0}}
(X^{(1)}(t-s)-X^{(2)}(t-s))^{2}ds]\text{.}%
\end{array}
\]

\end{description}
\end{lemma}

\noindent{Proof.} \quad By definition of the norms and standard properties of
the complex exponential function, we have%
\begin{align*}
&  ||\mathcal{L}(X^{(1)},X^{(2)})-\mathcal{L}(\widetilde{X}^{(1)}%
,\widetilde{X}^{(2)})||_{\mathcal{M}_{0}^{2}}^{2}\\
&  :={%
{\textstyle\int_{\mathcal{\mathbb{R}}^{2}}}
}|\widehat{\mathcal{L}}(X^{(1)},X^{(2)})(y_{1},y_{2})-\widehat{\mathcal{L}%
}(\widetilde{X}^{(1)},\widetilde{X}^{(2)})(y_{1},y_{2})|^{2}e^{-y_{1}%
^{2}-y_{2}^{2}}dy_{1}dy_{2}\\
&  ={%
{\textstyle\int_{\mathbb{R}^{2}}}
}|{%
{\textstyle\int_{\mathbb{R}^{2}}}
}e^{-i(x^{(1)}y_{1}+x^{(2)}y_{2})}d\mathcal{L}(X^{(1)},X^{(2)})(x^{(1)}%
,x^{(2)})\\
&  -%
{\textstyle\int_{\mathcal{\mathbb{R}}^{2}}}
e^{-i(\widetilde{x}^{(1)}y_{1}+\widetilde{x}^{(2)}y_{2})}d\mathcal{L}%
(\widetilde{X}^{(1)},\widetilde{X}^{(2)})(\widetilde{x}^{(1)},\widetilde
{x}^{(2)})|^{2}e^{-y_{1}^{2}-y_{2}^{2}}dy_{1}dy_{2}\\
&  ={%
{\textstyle\int_{\mathcal{\mathbb{R}}^{2}}}
}|\mathbb{E}[e^{-i(X^{(1)}y_{1}+X^{(2)}y_{2})}-e^{-i(\widetilde{X}^{(1)}%
y_{1}+\widetilde{X}^{(2)}y_{2})}]|^{2}e^{-y_{1}^{2}-y_{2}^{2}}dy_{1}dy_{2}\\
&  \leq{%
{\textstyle\int_{\mathcal{\mathbb{R}}^{2}}}
}\mathbb{E}[|e^{-i(X^{(1)}y_{1}+X^{(2)}y_{2})}-e^{-i(\widetilde{X}^{(1)}%
y_{1}+\widetilde{X}^{(2)}y_{2})}|^{2}]e^{-y_{1}^{2}-y_{2}^{2}}dy_{1}dy_{2}\\
&  ={%
{\textstyle\int_{\mathcal{\mathbb{R}}^{2}}}
}\mathbb{E}[(\cos(X^{(1)}y_{1}+X^{(2)}y_{2})-\cos(\widetilde{X}^{(1)}%
y_{1}+\widetilde{X}^{(2)}y_{2})^{2}\\
&  +(\sin(X^{(1)}y_{1}+X^{(2)}y_{2})-\sin(\widetilde{X}^{(1)}y_{1}%
+\widetilde{X}^{(2)}y_{2}))^{2}]e^{-y_{1}^{2}-y_{2}^{2}}dy_{1}dy_{2}\\
&  \leq{%
{\textstyle\int_{\mathcal{\mathbb{R}}^{2}}}
}(\mathbb{E}[|(X^{(1)}-\widetilde{X}^{(1)})y_{1}+(X^{(2)})-\widetilde{X}%
^{(2)})y_{2}|^{2}]\\
&  +\mathbb{E}[(X^{(1)}-\widetilde{X}^{(1)})y_{1}+(X^{(2)})-\widetilde
{X}^{(2)})y_{2}|^{2})]e^{-y_{1}^{2}-y_{2}^{2}}dy_{1}dy_{2}\\
&  =2{%
{\textstyle\int_{\mathcal{\mathbb{R}}^{2}}}
}(\mathbb{E}[|(X^{(1)}-\widetilde{X}^{(1)})y_{1}+(X^{(2)})-\widetilde{X}%
^{(2)})y_{2}|]^{2})e^{-y_{1}^{2}-y_{2}^{2}}dy_{1}dy_{2}\\
&  \leq4{%
{\textstyle\int_{\mathcal{\mathbb{R}}^{2}}}
}(\mathbb{E}[(X^{(1)}-\widetilde{X}^{(1)})^{2}]y_{1}^{2}+\mathbb{E}%
[(X^{(2)}-\widetilde{X}^{(2)})^{2}]y_{2}^{2})e^{-y_{1}^{2}-y_{2}^{2}}%
dy_{1}dy_{2}\\
&  \leq C_{0}\mathbb{E}[(X^{(1)}-\widetilde{X}^{(1)})^{2}+(X^{(2)}%
)-\widetilde{X}^{(2)})^{2}].
\end{align*}

\noindent Similarly, we get that
\[%
\begin{array}
[c]{lll}%
||\mathcal{L}(X_{t}^{(1)})-\mathcal{L}(X_{t}^{(2)})||_{\mathcal{M}_{0,\delta
}^{2}}^{2} & \leq &
{\textstyle\int_{\mathbb{-\delta}}^{0}}
\left\Vert \mathcal{L}(X^{(1)}(t-s))-\mathcal{L}(X^{(2)}(t-s))\right\Vert
_{\mathcal{M}_{0}^{2}}^{2}ds\\
& \leq & C_{0}\ \mathbb{E}[%
{\textstyle\int_{\mathbb{-\delta}}^{0}}
(X^{(1)}(t-s)-X^{(2)}(t-s))^{2}ds].
\end{array}
\]
$\square$\newline

\subsection{Spaces}

\noindent Throughout this work, we will use the following spaces:

\begin{itemize}
\item $\mathbb{L}^{2}$ is the space of measurable functions{ }$\sigma
:[0,\delta]\rightarrow\mathbb{R}$, such that
\[
\parallel\sigma\parallel_{\mathbb{L}^{2}}^{2}:=%
{\textstyle\int_{0}^{\delta}}
|\sigma(r)|^{2}dr<\infty.
\]

\item $\mathcal{S}^{2}$ is the set of ${\mathbb{R}}$-valued $\mathbb{F}%
$-adapted c\`{a}dl\`{a}g processes $(X(t))_{t\in\lbrack0,T]}$, such that
\[
{\Vert X\Vert}_{\mathcal{S}^{2}}^{2}:={\mathbb{E}}[\sup_{t\in\lbrack
0,T]}|X(t)|^{2}]~<~\infty\;.
\]

\item $L^{2}$ is the set of ${\mathbb{R}}$-valued $\mathbb{F}$-predictable
processes $(Q(t))_{t\in\lbrack0,T]}$, such that
\[
\Vert Q\Vert_{L^{2}}^{2}:={\mathbb{E}}[%
{\textstyle\int_{0}^{T}}
|Q(t)|^{2}dt]<~\infty\;.
\]

\item $\Xi$ is the set of $\mathbb{G}$-adapted, nondecreasing
right-continuous processes $\xi$ with $\xi(0^{-})=0$ (the set of admissible
singular controls).

\item $L^{2}(\Omega,\mathcal{F}_{t})$ is the set of ${\mathbb{R}}$-valued
square integrable $\mathcal{F}_{t}$-measurable random variables.

\item $\mathcal{R}$ is the set of functions $r:\mathbb{R}%
_{0}\rightarrow\mathbb{R}.$

\item $C_{a}([0,T],\mathcal{M}_{0})$ denotes the set of absolutely continuous
functions $m:[0,T]\rightarrow\mathcal{M}_{0}.$

\end{itemize}

\section{Existence and Uniqueness of Solutions of Reflected AMBSDEs}

In this section, we will prove existence and uniqueness of solutions of
reflected mean-field BSDEs with a generator which is \emph{(time-) advanced},
in the sense that at any time $t$, the generator may depend on future values
up to a positive constant $\delta$ as follows:\newline For a given driver $F$,
terminal value $R$ and barrier (or obstacle) process $S$, we say that an
$\mathbb{F}$-adapted process $(Y,Z,K)\in\mathcal{S}^{2}\times L^{2}\times\Xi$
is a solution of the corresponding reflected AMBSDEs if the following holds:%
\begin{equation}%
\left\{
\begin{array}
[c]{l}%
(i)Y(t)=R+%
{\textstyle\int_{t}^{T}}
F(s,Y(s),Z(s),\mathbb{E}[Y^{s}|\mathcal{F}_{s}],\mathbb{E}[Z^{s}%
|\mathcal{F}_{s}],\mathcal{L}(Y^{s},Z^{s}))ds\\
\quad\quad\quad\quad\quad\quad+K(T)-K(t)-%
{\textstyle\int_{t}^{T}}
Z(s)dB(s);\quad0\leq t\leq T,\\
(ii)Y(t)\geq S(t);\quad0\leq t\leq T,\\
(iii){%
{\textstyle\int_{0}^{T}}
}(Y(t)-S(t))dK^{c}(t)=0\text{ a.s. and }\triangle K^{d}(t)=-\triangle
Y(t)\mathbf{1}_{\{Y(t^{-})=S(t^{-})\}}\text{ a.s.},\\
(iv)Y(t)=R;\quad t\geq T,\\
(v)Z(t)=0;\quad t>T,
\end{array}
\right.
\label{eq3.1}%
\end{equation}
where $Y^{s}=\left(  Y(s+r)\right)  _{r\in\left[  0,\delta\right]  }%
,Z^{s}=\left(  Z(s+r)\right)  _{r\in\left[  0,\delta\right]  },$ the terminal
condition $R\in L^{2}(\Omega,\mathcal{F}_{T})$, the driver $F:[0,T]\times
\Omega\times\mathbb{R}^{2}\mathbb{\times L}^{2}\times\mathbb{L}^{2}%
\times\mathcal{M}_{0,\delta}\longrightarrow\mathbb{R}$ is $\mathcal{F}_{t}%
$-progressively measurable and we have denoted by $K^{c}$ and $K^{d}$ the
continuous and discontinuous parts of $K$ respectively.

\noindent We may remark here that in order to guarantee adaptedness, the
time-advanced terms are given under conditional expectation with respect to
$\mathcal{F}_{s}$.

\noindent Our result can be regarded as an extension of the existing results
on advanced BSDEs of Peng \& Yang \cite{peng2009}, \O ksendal \textit{et al}
\cite{Oksendal2011}, Jeanblanc \textit{et al} \cite{JLA} and we refer here to
the paper by Quenez and Sulem \cite{QS} on reflected BSDEs for c\`{a}dl\`{a}g obstacle.

\noindent To obtain the existence and the uniqueness of a solution, we make
the following set of assumptions:

\begin{itemize}
\item For the driver $F,$ we assume

\item[ (i)] There exists a constant $c\in%
\mathbb{R}
$ such that
\[
|F(\cdot,0,0,0,0,\mathcal{L}(0,0))|\leq c,
\]
where $\mathcal{L}(0,0)$ is the Dirac measure with mass at zero.

\item[ (ii)] There exists a constant $C_{Lip}^{F}\in%
\mathbb{R}
$ such that, for $t\in\lbrack0,T],$
\begin{align*}
&  |F(t,y_{1},z_{1},y_{2},z_{2},\mathcal{L}(y_{2},z_{2}))-F(t,y_{1}^{\prime
},z_{1}^{\prime},y_{2}^{\prime},z_{2}^{\prime},\mathcal{L}(y_{2}^{\prime
},z_{2}^{\prime}))|^{2}\\
&  \leq C_{Lip}^{F}\{|y_{1}-y_{1}^{\prime}|^{2}+|z_{1}-z_{1}^{\prime}%
|^{2}+||y_{2}-y_{2}^{\prime}||_{\mathbb{L}^{2}}^{2}+||z_{2}-z_{2}^{\prime
}||_{\mathbb{L}^{2}}^{2}\\
&  +||\mathcal{L}(y_{2},z_{2})-\mathcal{L}(y_{2}^{\prime},z_{2}^{\prime
})||_{\mathcal{M}_{0,\delta}}^{2})\},
\end{align*}
for all $y_{1},z_{1},y_{1}^{\prime},z_{1}^{\prime}\in{\mathbb{R}},$
$y_{2},z_{2},y_{2}^{\prime},z_{2}^{\prime}\in\mathbb{L}^{2},$ $\mathcal{L}%
(y_{2},z_{2}),\mathcal{L}(y_{2}^{\prime},z_{2}^{\prime})\in\mathcal{M}%
_{0,\delta}.$

\item For the barrier $S,$ we assume:

\item[(iii)] The barrier $S$ is nondecreasing, $\mathbb{F}$-adapted,
c\`{a}dl\`{a}g process satisfying
\[
\mathbb{E}[\underset{t\in\lbrack0,T]}{\sup}|S(t)|^{2}]<\infty.
\]

\item[(iv)] $Y(t)\geq S(t); 0\leq t\leq T$.

\item For the local time $K,$ we assume:

\item[(v)] $K$ is a nondecreasing ${\mathbb{F}}$-adapted c\`{a}dl\`{a}g
process with $K(0^{-})=0,$ such that ${%
{\textstyle\int_{0}^{T}}
}(Y(t)-S(t))dK^{c}(t)=0$ a.s. and $\triangle K^{d}(t)=-\triangle
Y(t)\mathbf{1}_{\{Y(t^{-})=S(t^{-})\}}$ a.s.
\end{itemize}

\begin{theorem}
[Existence and Uniqueness]Under the above assumptions (i)-(v), the reflected
AMBSDEs (\ref{eq3.1}) has a unique solution $(Y,Z,K) \in\mathcal{S}^{2} \times
L^{2} \times\Xi$.
\end{theorem}

\noindent{Proof.} \quad For $t\in\lbrack0,T]$ and for all {$\beta>0$, we
define the Hilbert space }$\mathbb{H}_{{\beta}}^{2}$ to be the set of all
$(Y,Z)\in\mathcal{S}^{2}\times L^{2}$, equipped with the norm
\[
||(Y,Z)||_{\mathbb{H}_{{\beta}}^{2}}^{2}:={\mathbb{E}}[%
{\textstyle\int_{0}^{T+\delta}}
e^{\beta t}(Y^{2}(t)+Z^{2}(t))dt]\;.
\]

\noindent Define the {mapping $\Phi:\mathbb{H}_{{\beta}}^{2}\rightarrow$%
}$\mathbb{H}_{{\beta}}^{2}$ by $\Phi(y,z)=(Y,Z)$ where $(Y,Z)\in
${{$\mathcal{S}^{2}$}$\times${$L^{2}(\subset{{L^{2}}\times{L^{2}}})$}} is
defined by%
\[
\left\{
\begin{array}
[c]{ll}%
Y(t) & =R+{%
{\textstyle\int_{t}^{T}}
}F(s,{y(s),z(s)},\mathbb{E}[y^{s}|\mathcal{F}_{s}],\mathbb{E}[z^{s}%
|\mathcal{F}_{s}],\mathcal{L}(y^{s},z^{s}))ds\\
& +K(T)-K(t)-{%
{\textstyle\int_{t}^{T}}
}Z(s)dB(s); \quad0\leq t\leq T,\\
Y(t) & =R; \quad t\geq T,\\
Z(t) & =0; \quad t>T.
\end{array}
\right.
\]
To prove the theorem, it suffices to prove that {$\Phi$ is a contraction
mapping in }$\mathbb{H}_{{\beta}}^{2}${ under the norm $||\cdot||_{\mathbb{H}%
_{{\beta}}^{2}}$ for large enough }${\beta}${. For two arbitrary elements
$(y_{1},z_{1},k_{1})$ and $(y_{2},z_{2},k_{2})$, we denote their difference
by}
\[
(\widetilde{y},\widetilde{z},\widetilde{k})=(y_{1}-y_{2},z_{1}-z_{2}%
,k_{1}-,k_{2})\;.
\]
Applying It\^{o} formula for semimartingale, we get%
\begin{align*}
&  \mathbb{E}[%
{\textstyle\int_{0}^{T}}
e^{\beta t}(\beta\widetilde{Y}^{2}(t)+\widetilde{Z}^{2}(t))dt]\\
&  =2\mathbb{E}[%
{\textstyle\int_{0}^{T}}
e^{\beta t}\widetilde{Y}(t)\{F(t,{y_{1}(t),z_{1}(t)},\mathbb{E}[y_{1}%
^{t}|\mathcal{F}_{t}],\mathbb{E}[z_{1}^{t}|\mathcal{F}_{t}],\mathcal{L}%
(y_{1}^{t},z_{1}^{t}))\\
&  -F(t,{y_{2}(t),z_{2}(t)},\mathbb{E}[y_{2}^{t}|\mathcal{F}_{t}%
],\mathbb{E}[z_{2}^{t}|\mathcal{F}_{t}],\mathcal{L}(y_{2}^{t},z_{2}%
^{t}))\}dt]\;\\
&  +2\mathbb{E}[%
{\textstyle\int_{0}^{T}}
e^{\beta t}\widetilde{Y}(t)dK^{1}(t)]-2\mathbb{E}[%
{\textstyle\int_{0}^{T}}
e^{\beta t}\widetilde{Y}(t)dK^{2}(t)].
\end{align*}
We have that
\begin{align*}
\widetilde{Y}(t)dK^{1,c}(t)  & =(Y^{1}(t)-S(t))dK^{1,c}(t)-(Y^{2}%
(t)-S(t))dK^{1,c}(t)\\
& =-(Y^{2}(t)-S(t))dK^{1,c}(t)\leq0\text{ a.s.,}%
\end{align*}
and by symmetry, we have also $\widetilde{Y}(t)dK^{2,c}(t)\geq0$ a.s. For the
discontinuous case, we have as well%
\begin{align*}
\widetilde{Y}(t)dK^{1,d}(t)  & =(Y^{1}(t)-S(t))dK^{1,d}(t)-(Y^{2}%
(t)-S(t))dK^{1,d}(t)\\
& =-(Y^{2}(t)-S(t))dK^{1,d}(t)\leq0\text{ a.s.,}%
\end{align*}
and by symmetry, we have also $\widetilde{Y}(t)dK^{2,d}(t)\geq0$ a.s.

\noindent{By Lipschitz assumption and standard estimates, it follows that }%

\begin{align*}
&\mathbb{E}[
{\textstyle\int_{0}^{T}}
e^{\beta t}(\beta\widetilde{Y}^{2}(t)+\widetilde{Z}^{2}(t))dt]\\
& \leq 8\rho C^{2}\text{ }\mathbb{E}[
{\textstyle\int_{0}^{T}}
e^{\beta t}\widetilde{Y}^{2}(t)dt]\\
&+\tfrac{1}{2\rho}\mathbb{E}[{\textstyle\int_{0}^{T}}
e^{\beta t}(\widetilde{y}^{2}(t)+\widetilde{z}^{2}(t)+{\textstyle\int_{0}^{\delta}}
(\widetilde{y}^{2}(t+r)+\widetilde{z}^{2}(t+r))dr)dt]\;.
\end{align*}
By change of variable $s=t+r$, we get
\begin{align*}
& \mathbb{E}[
{\textstyle\int_{0}^{T}}
e^{\beta t}
{\textstyle\int_{0}^{\delta}}
(\widetilde{y}^{2}(t+r)+\widetilde{z}^{2}(t+r))dr)dt]\\
& \leq\mathbb{E}[
{\textstyle\int_{0}^{T}}
e^{\beta t}
{\textstyle\int_{t}^{t+\delta}}
(\widetilde{y}^{2}(s)+\widetilde{z}^{2}(s))ds)dt].
\end{align*}
Fubini's theorem gives that
\begin{align*}
& \mathbb{E}[%
{\textstyle\int_{0}^{T}}
e^{\beta t}%
{\textstyle\int_{0}^{\delta}}
(\widetilde{y}^{2}(t+r)+\widetilde{z}^{2}(t+r))dr)dt]\\
& \leq\mathbb{E}[%
{\textstyle\int_{0}^{T+\delta}}
(%
{\textstyle\int_{s-\delta}^{s}}
e^{\beta t}dt)(\widetilde{y}^{2}(s)+\widetilde{z}^{2}(s)))ds]\\
& \leq\mathbb{E}[%
{\textstyle\int_{0}^{T+\delta}}
e^{\beta s}(\widetilde{y}^{2}(s)+\widetilde{z}^{2}(s)))ds].
\end{align*}
Consequently, by choosing $\beta=1+8\rho C^{2},$ we have
\[
\mathbb{E}[%
{\textstyle\int_{0}^{T}}
e^{\beta t}(\widetilde{Y}^{2}(t)+\widetilde{Z}^{2}(t))dt]\leq\tfrac{1}{\rho
}\;\mathbb{E}[%
{\textstyle\int_{0}^{T+\delta}}
e^{\beta t}(\widetilde{y}^{2}(t)+\widetilde{z}^{2}(t))dt]\;.
\]
Since $\widetilde{Y}(t)=\widetilde{Z}(t)=0$ for $t>T$, we get
\[
||(\widetilde{Y},\widetilde{Z})||_{\mathbb{H}_{{\beta}}^{2}}^{2}\leq\tfrac
{1}{\rho}\;||(\widetilde{y},\widetilde{z})||_{\mathbb{H}_{{\beta}}^{2}}^{2}\;.
\]
{ }For{ }$\rho${$>1$, we get that }$\Phi$ is a contraction on $\mathbb{H}%
_{{\beta}}^{2}.$
\fproof
\section{Reflected AMBSDEs and optimal stopping under partial information}

In this section we recall a connection between reflected AMBSDEs and optimal
stopping problems under partial information.

\begin{definition}
Let $F:\Omega\times\lbrack0,T]\times\mathbb{R}^{2}\times
\mathbb{L}^{2}\times\mathbb{L}^{2}\times\mathcal{M}_{0,\delta}\rightarrow\mathbb{R}$ be a given function.
 \\Assume that:

$\bullet$ $F$ is $\mathbb{G}$-adapted and $|F(t,0,0,0,0,\mathcal{L}(0,0))|<c$,  for all $t$; for some constant $c$.

$\bullet$ $S(t)$ is a given $\mathbb{F}$-adapted c\`{a}dl\`{a}g nondecreasing
process, such that \newline%
\[
\mathbb{E[}\underset{t\in\lbrack0,T]}{\sup}(S(t))^{2}]<\infty.
\]

$\bullet$ The terminal value $R \in L^{2}\left(  \Omega,\mathcal{F}%
_{T}\right)  $ is such that $R \geq S(T)$ a.s.

We say that a $\mathbb{G}$-adapted triplet $\left(  Y,\mathcal{Z},K\right)  $ is a
solution of the reflected AMBSDE with driver $F$, terminal value $R$ and the
reflecting barrier $S(t)$ under the filtration $\mathbb{G}$, if the following hold:

\begin{enumerate}
\item
\[
\mathbb{E[}%
{\textstyle\int_{0}^{T}}
|F(s,Y(s),Z(s),\mathbb{E[}Y^{s}|\mathcal{F}_{s}],\mathbb{E[}\mathcal{Z}^{s}%
|\mathcal{F}_{s}],\mathcal{L}(Y^{s},\mathcal{Z}^{{s}}))|^{2}ds]<\infty,
\]

\item
\[
\mathcal{Z}(t) \text{  is a } \mathbb{G}-martingale,
\]

\item
\[%
\begin{array}
[c]{c}%
Y(t)=R+%
{\textstyle\int_{t}^{T}}
F(s,Y(s),Z(s),\mathbb{E[}Y^{s}|\mathcal{F}_{s}],\mathbb{E[}\mathcal{Z}^{s}%
|\mathcal{F}_{s}],\mathcal{L}(Y^{s},\mathcal{Z}^{{s}}))ds\\
-%
{\textstyle\int_{t}^{T}}
dK(s)-%
{\textstyle\int_{t}^{T}}
d\mathcal{Z}(s); \quad t\in\left[  0,T\right]  ,
\end{array}
\newline%
\]
or, equivalently,\newline%
\[%
\begin{array}
[c]{c}%
Y(t)=\mathbb{E}[R+%
{\textstyle\int_{t}^{T}}
F(s,Y(s),Z(s),\mathbb{E[}Y^{s}|\mathcal{F}_{s}],\mathbb{E[}\mathcal{Z}^{s}%
|\mathcal{F}_{s}],\mathcal{L}(Y^{s},\mathcal{Z}^{{s}}))ds\\
-%
{\textstyle\int_{t}^{T}}
dK(s)|\mathcal{G}_{t}]; t\in\left[  0,T\right]  ,
\end{array}
\]

\item $K(t)$ is nondecreasing, $\mathbb{G}$-adapted, c\`{a}dl\`{a}g and $K(0^{-})=0,$

\item $Y(t)\geq S(t)$ a.s.; $t\in\lbrack0,T],$

\item $%
{\textstyle\int_{0}^{T}}
(Y(t)-S(t))dK(t)=0$ a.s.
\end{enumerate}
\end{definition}

\vskip0.3cm The following result is essentially due to El Karoui \textit{et
al} \cite{EKPPQ}. See also \O ksendal \& Sulem \cite{OS3} and \O ksendal \& Zhang \cite{OZ}.

\begin{theorem}
For $t\in\lbrack0,T]$, let $\mathcal{T}_{[t,T]}$ denote the set of all
$\mathbb{G}$-stopping times $\tau:\Omega\mapsto\lbrack t,T].$\newline Suppose
$\left(  Y,\mathcal{Z},K\right)  $ is a solution of the reflected AMBSDE above.
\end{theorem}


\begin{description}
\item[(i)] Then $Y(t)$ is the solution of the optimal stopping problem%
\[%
\begin{array}
[c]{c}%
Y(t)=\underset{\tau\in\mathcal{T}_{[t,T]}}{ess\sup}\quad\{\mathbb{E}[%
{\int_{t}^{\tau}}
F(s,Y(s),\mathcal{Z}(s),Y^{s},\mathcal{Z}^{s},\mathcal{L}(Y^{s},\mathcal{Z}^{s}))ds\\
+S(\tau)\mathbf{1}_{\tau<T}+R\mathbf{1}_{\tau=T}|\mathcal{G}_{t}]\}; \quad
t\in\left[  0,T\right]  .
\end{array}
\]

\item[(ii)] Moreover the solution process $K(t)$ is given by%
\begin{align}
K(T)-K(T-t)&=\underset{s\leq t}{\max}\Big\{R+%
\int_{T-s}^{T}
F(r,Y(r),\mathcal{Z}(r),\mathbb{E[}Y^{r}|\mathcal{F}_{r}],\mathbb{E[}\mathcal{Z}^{r}|\mathcal{F}_{r}],\mathcal{L}(Y^{r},\mathcal{Z}^{r}))dr\nonumber\\
&-\int _{T-s}^{T} d\mathcal{Z}(r)-S(T-s)\Big\}^{-}; \quad t\in\left[
0,T\right]  ,
\end{align}
where $x^{-}=\max(-x,0),$
 and an optimal stopping time $\hat{\tau}_{t}$ is given by%
\begin{align*}
\hat{\tau}_{t}:&=\inf\{s\in\lbrack t,T],Y(s)\leq S(s)\}\wedge T\\
&  =\inf\{s\in\lbrack t,T],K(s) > K(t)\}\wedge T.
\end{align*}
\item[(iii)] In particular, if we choose $t=0$, we get that
\begin{align*}
\hat{\tau}_{0}:&=\inf\{s\in\lbrack 0,T],Y(s)\leq S(s)\}\wedge T\\
&  =\inf\{s\in\lbrack 0,T],K(s) > 0\}\wedge T,
\end{align*}
solves the optimal stopping problem

\begin{align*}
Y(0)&=\sup_{\tau\in\mathcal{T}_{[0,T]}}\mathbb{E}[{\textstyle\int_{0}^{\tau}}
F(s,Y(s), Z(s), \mathbb{E}[Y^{s}|\mathcal{F}_s],\mathbb{E}[Z^{s}|\mathcal{F}_s],\mathcal{L}(Y^{s},Z^{s}))ds\\
&+S(\tau)\mathbf{1}_{\tau<T}%
+R\mathbf{1}_{\tau=T}]; t\in\left[  0,T\right]  .
\end{align*}

\end{description}

\section{Optimal singular control of memory mean-field SDEs}

We now return to the singular control problem stated in the
Introduction:\newline

\subsection{Problem statement}

Consider the following \emph{mean-field memory singular controlled system},
with a state process $X(t)=X^{\xi}(t)$ and a singular control process
$\xi(t),$ of the form
\begin{equation}
\left\{
\begin{array}
[c]{l}%
dX(t)=b(t,X(t),X_{t},M(t),M_{t},\xi(t))dt+\sigma(t,X(t),X_{t},M(t),M_{t}%
,\xi(t))dB(t)\\
\quad\quad\quad+\lambda(t)d\xi(t); \quad t\in\lbrack0,T],\\
X(t)=\alpha(t); \quad t\in\lbrack-\delta,0],
\end{array}
\right.  \label{eq6.1}%
\end{equation}
where $X_{t}=\{X(t-s)\}_{0\leq s\leq\delta},$ $M(t)=\mathcal{L}(X(t)),$
$M_{t}=\{M(t-s)\}_{0\leq s\leq\delta},$ $b,\sigma:\Omega\times\lbrack
0,T]\times%
\mathbb{R}
\times\mathbb{L}^{2}\times\mathcal{M}_{0}\times\mathcal{M}_{0,\delta}%
\times{\mathbb{R}\times\Xi\rightarrow\mathbb{R}},$ $\lambda:[0,T]\rightarrow%
\mathbb{R}
.$\newline We assume that our control process $\xi(t)$ is ${\mathbb{R}}%
$-valued right-continuous $\mathbb{G}$-adapted processes, and $t\mapsto\xi(t)$
is increasing (nondecreasing) with $\xi(0^{-})=0$, and such that the
corresponding state equation has a unique solution $X$ with $\omega\mapsto
X(t,\omega)\in L^{2}(\mathbb{P})$ for all $t$. The set of such processes $\xi$
is denoted by $\Xi$. \newline

\noindent The \emph{performance functional} is assumed to be of the form
\begin{equation}%
\begin{array}
[c]{c}%
J(\xi)={\mathbb{E}}[%
{\textstyle\int_{0}^{T}}
f(t,X(t),X_{t},M(t),M_{t},\xi(t))dt+g(X(T),M(T))\\
\qquad+%
{\textstyle\int_{0}^{T}}
h(t,X(t))d\xi(t)];\quad\xi\in\Xi,
\end{array}
\label{eq6.3}%
\end{equation}
where $f:\Omega\times\lbrack0,T]\times%
\mathbb{R}
\times\mathbb{L}^{2}\times\mathcal{M}_{0}\times\mathcal{M}_{0,\delta}%
\times{\mathbb{R}\times\Xi\rightarrow\mathbb{R}},$ $h:\Omega\times
\lbrack0,T]\times%
\mathbb{R}
\rightarrow%
\mathbb{R}
,$ $g:\Omega\times%
\mathbb{R}
\times\mathcal{M}_{0}\rightarrow%
\mathbb{R}
.$\newline The problem we consider, is the following:

\begin{problem}
\label{prob} Find an optimal control $\hat{\xi}\in\Xi$, such that
\begin{equation}
J(\hat{\xi})=\sup_{\xi\in\Xi}J(\xi)\,. \label{eq6.4}%
\end{equation}

\end{problem}

\noindent First we explain some notation and introduce some useful dual
operators. \newline Let $L_{0}^{2}$ denote the set of measurable stochastic
processes $Y(t)$ on $\mathbb{R}$ such that $Y(t)=0$ for $t<0$ and for $t>T$
and
\[
\EE[{\textstyle\int_{0}^{T}} Y^{2}(t)dt]<\infty\quad a.s.
\]

\begin{itemize}
\item Let $G(t,\bar{x})=G_{\bar{x}}(t,\cdot):[0,T]\times\mathbb{L}^{2}%
\mapsto\mathbb{R}$ be a bounded linear functional on $\mathbb{L}^{2}$ for each
$t$, uniformly bounded in $t$.Then the map
\[
Y\mapsto\mathbb{E}[%
{\textstyle\int_{0}^{T}}
\left\langle G_{\overline{x}}(t),Y_{t}\right\rangle dt];\quad Y\in L_{0}^{2}%
\]
is a bounded linear functional on $L_{0}^{2}$. Therefore, by the Riesz
representation theorem there exists a unique process denoted by $G_{\bar{x}%
}^{\ast}(t)\in L_{0}^{2}$ such that
\begin{equation}
{\mathbb{E}}[{%
{\textstyle\int_{0}^{T}}
}\left\langle G_{\overline{x}}(t),Y_{t}\right\rangle dt]={\mathbb{E}}[{%
{\textstyle\int_{0}^{T}}
}G_{\bar{x}}^{\ast}(t)Y(t)dt], \label{eq6.7a}%
\end{equation}
for all $Y\in L_{0}^{2}$.
\end{itemize}

We illustrate these operators by some auxiliary results.

\begin{lemma}
Consider the case when
\[
G_{\bar{x}}(t,\cdot)=\left\langle F,\cdot\right\rangle p(t), \text{ with } p \in L_0^2.
\]
Then
\begin{equation}
G_{\bar{x}}^{\ast}(t):=\left\langle F,p^{t}\right\rangle \label{eq4.8}%
\end{equation}
satisfies \eqref{eq6.7a}, where $p^{t}:=\{p(t+r)\}_{r\in\lbrack0,\delta]}$.
\end{lemma}

\noindent{Proof.} \quad\quad We must verify that if we define $G^{*}_{\bar{x}%
}(t)$ by \eqref{eq4.8}, then \eqref{eq6.7a} holds. To this end, choose $Y\in
L_{0}^{2}$ and consider%

\begin{align}
&
{\textstyle\int_{0}^{T}}
\left\langle F,p^{t}\right\rangle Y(t)dt={%
{\textstyle\int_{0}^{T}}
}\left\langle F,\{p(t+r)\}_{r\in\lbrack0,\delta]}\right\rangle
Y(t)dt\nonumber\\
&  =%
{\textstyle\int_{0}^{T}}
\left\langle F,\{Y(t)p(t+r)\}_{r\in\lbrack0,\delta]}\right\rangle
dt=\left\langle F,\{%
{\textstyle\int_{r}^{T+r}}
Y(u-r)p(u)du\}_{r\in\lbrack0,\delta]}\right\rangle \nonumber\\
&  =\left\langle F,\{%
{\textstyle\int_{0}^{T}}
Y(u-r)p(u)du\}_{r\in\lbrack0,\delta]}\right\rangle ={%
{\textstyle\int_{0}^{T}}
}\left\langle F,Y_{u}\right\rangle p(u)du\nonumber\\
&  ={%
{\textstyle\int_{0}^{T}}
}\left\langle \nabla_{\bar{x}}G(u),Y_{u}\right\rangle du.\nonumber
\end{align}

\hfill$\square$ \bigskip

\begin{example}
(i) For example, if $a\in\mathbb{R}^{[0,\delta]}$ is a bounded function and
$F(\bar{x})$ is the averaging operator defined by
\[
F(\bar{x})=\left\langle F,\bar{x}\right\rangle =%
{\textstyle\int_{-\delta}^{0}}
a(s)x(s)ds
\]
when $\bar{x}=\{x(s)\}_{s\in\lbrack0,\delta]}$, then
\[
\left\langle F,p^{t}\right\rangle ={%
{\textstyle\int_{0}^{\delta}}
}a(r)p(t+r)dr.
\]
(ii) Similarly, if $t_{0}\in\lbrack0,\delta]$ and $G$ is evaluation at $t_{0}%
$, i.e.
\[
G(\bar{x})=x(t_{0})\text{ when }\bar{x}=\{x(s)\}_{s\in\lbrack0,\delta]},
\]
then
\[
\left\langle G,p^{t}\right\rangle =p(t+t_{0}).
\]

\end{example}

We now have the machinery to start working on Problem (\ref{prob}%
)\textbf{.}\newline Let $\widehat{\mathcal{M}}$ be the set of all random
measures on $[0,T]$. Define the \textit{(singular) Hamiltonian}
\[
H:[0,T]\times\mathbb{R}\times\mathbb{L}^{2}\times\mathcal{M}_{0}%
\times\mathcal{M}_{0,\delta}\times\Xi\times\mathbb{R}\times\mathbb{R}\times
C_{a}([0,T],\mathcal{M}_{0})\mapsto\widehat{\mathcal{M}}%
\]
as the following random measure:
\begin{align}
dH(t)  &  =dH(t,x,\bar{x},m,\bar{m},\xi,p^{0},q^{0},p^{1})\label{eq5.2a}\\
&  =H_{0}(t,x,\bar{x},m,\bar{m},\xi,p^{0},q^{0},p^{1})dt+\{\lambda
(t)p^{0}+h(t,x)\}d\xi(t)\,,\nonumber
\end{align}
where
\begin{align}
&  H_{0}(t,x,\bar{x},m,\bar{m},\xi,p^{0},q^{0},p^{1})\label{eq5.3a}\\
&  :=f(t,x,\bar{x},m,\bar{m},\xi)+b(t,x,\bar{x},m,\bar{m},\xi)p^{0}%
+\sigma(t,x,\bar{x},m,\bar{m},\xi)q^{0}+\left\langle p^{1},\beta
(m)\right\rangle ,\nonumber
\end{align}
where $\beta(m)$ is defined below. Here $m$ denotes a generic value of the
measure $M(t)$.
We assume that f, b, $\sigma,\gamma,h$ and $g$ are Fr\'{e}chet differentiable
$(C^{1})$ in the variables $x,\bar{x},m,\bar{m},\xi$. Then the same holds for
$H_{0}$ and $H$.\newline\noindent We define the adjoint processes
$(p^{0},q^{0}),(p^{1},q^{1})$ as the solutions of the following BSDEs, respectively:%

\begin{equation}
\left\{
\begin{array}
[c]{ll}%
dp^{0}(t) & =-\Big\{\tfrac{\partial H_{0}}{\partial x}(t)+\mathbb{E[}%
\nabla_{\bar{x}}^{\ast}H_{0}(t)|\mathcal{F}_{t}]\Big\}dt-\frac{\partial
h}{\partial x}(t)d\xi(t)+q^{0}(t)dB(t);\quad t\in\lbrack0,T],\\
p^{0}(t) & =\tfrac{\partial g}{\partial x}(T);\quad t\geq T,\\
q^{0}(t) & =0;\quad t>T,
\end{array}
\right.  \label{eqp0}%
\end{equation}

and
\begin{equation}
\left\{
\begin{array}
[c]{ll}%
dp^{1}(t) & =-\{\nabla_{m}H_{0}(t)+\mathbb{E}\left[  \nabla_{\bar{m}}^{\ast
}H_{0}(t)|\mathcal{F}_{t}\right]  \}dt+q^{1}(t)dB(t);\quad t\in\lbrack0,T],\\
p^{1}(t) & =\nabla_{m}g(T);\quad t\geq T,\\
q^{1}(t) & =0;\quad t>T,
\end{array}
\right.  \label{eqp1}%
\end{equation}
where $g(T)=g(X(T),M(T))$ and
\[
H_{0}(t)=H_{0}(t,x,\bar{x},m,\bar{m},\xi,p^{0},q^{0},p^{1})_{x=X(t),\bar
{x}=X_{t},m=M(t),\bar{m}=M_{t},\xi=\xi(t),p^{0}=p^{0}(t),q^{0}=q^{0}%
(t),p^{1}=p^{1}(t)}.
\]
Here $\nabla_{m}H_{0}$ is the Frech\'{e}t derivative of $H_{0}$ with respect
to $m$, and $\nabla_{\bar{m}}^{\ast}H_{0}$ is defined similarly to
$\nabla_{\bar{x}}^{\ast}H_{0}$.\newline



\subsection{A sufficient maximum principle for singular mean field control
with partial information}

\label{sec6.2}

We proceed to state a sufficient maximum principle (a verification theorem)
for the singular mean-field control problem described by \eqref{eq6.1} -
\eqref{eq6.4}. Because of the mean-field terms, it is natural to consider the
two-dimensional system $(X(t),M(t))$, where the dynamics for $M(t)$ is the
following:
\begin{equation}%
\begin{cases}
dM(t)=\beta(M(t)dt,\nonumber\\
M(0)\in\mathcal{M}_{0},
\end{cases}
\end{equation}
where we have put $\beta(M(t))=M^{\prime}(t)$. See Lemma \ref{m'}.\newline

\begin{theorem}
[\emph{Sufficient maximum principle for mean-field singular control}%
]\label{th5.1a} Let $\hat{\xi}\in\Xi$ be such that the system of \eqref{eq6.1}
and \eqref{eqp0} - \eqref{eqp1} has a solution $\hat{X}(t),\hat{p}^{0}%
(t),\hat{q}^{0}(t),\hat{p}^{1}(t),\hat{q}^{1}(t)$ and set $\hat{M}%
(t)=\mathcal{L}(\hat{X}(t))$.
Suppose the following conditions hold:

\begin{itemize}
\item (The concavity assumptions) The functions
\begin{align}
&  \mathbb{R}\times\mathbb{L}^{2}\times\mathcal{M}_{0}\times\mathcal{M}%
_{0,\delta}\times\Xi\ni(x,\bar{x},m,\bar{m},\xi)\mapsto dH(t,x,\bar{x}%
,m,\bar{m},\xi,\hat{p}^{0}(t),\hat{q}^{0}(t),\hat{p}^{1}(t),\hat{q}%
^{1}(t))\nonumber\\
&  \text{ and }\nonumber\\
&  \mathbb{R}\times\mathcal{M}_{0}\ni(x,m)\mapsto g(x,m)\nonumber\\
&  \text{ are concave for all }t\in\lbrack0,T]\text{ and almost all }\omega
\in\Omega. \label{eq3.10a}%
\end{align}

\item (Conditional variational inequality) For all $\xi\in\Xi$ we have
\[
\mathbb{E}[dH(t)|\mathcal{G}_{t}]\leq\mathbb{E}[d\hat{H}(t)|\mathcal{G}%
_{t}],\,
\]
i.e.
\begin{equation}%
\begin{array}
[c]{c}%
\mathbb{E}[H_{0}(t)|\mathcal{G}_{t}]dt+\mathbb{E}[\lambda(t)\hat{p}%
^{0}(t)+\hat{h}(t)|\mathcal{G}_{t}]d\xi(t)\\
\leq\mathbb{E}[\hat{H}_{0}(t)|\mathcal{G}_{t}]dt+\mathbb{E}[\lambda(t)\hat
{p}^{0}(t)+\hat{h}(t)|\mathcal{G}_{t}]d\hat{\xi}(t),
\end{array}
\label{eq5.21}%
\end{equation}
where the inequality is interpreted in the sense of inequality between random
measures in $\mathcal{M}$.
\end{itemize}

Then $\hat{\xi}(t)$ is an optimal control for $J(\xi)$.
\end{theorem}

\noindent{Proof.} \quad Choose $\xi\in\Xi$ and consider
\[
J(\xi)-J(\hat{\xi})=I_{1}+I_{2}+I_{3},
\]
where
\begin{align}
&  I_{1}={\mathbb{E}}[%
{\textstyle\int_{0}^{T}}
\{f(t)-\hat{f}(t)\}dt],\nonumber\\
&  I_{2}={\mathbb{E}}[g(T)-\hat{g}(T)],\nonumber\\
&  I_{3}={\mathbb{E}}[%
{\textstyle\int_{0}^{T}}
h(t)d\xi(t)-\hat{h}(t)d\hat{\xi}(t)].\label{eq6.16}%
\end{align}
By the definition of the Hamiltonian \eqref{eq5.3a} we have%
\begin{equation}%
\begin{array}
[c]{ll}%
I_{1} & =\mathbb{E[}%
{\textstyle\int_{0}^{T}}
\{H_{0}(t)-\hat{H}_{0}(t)-\hat{p}^{0}(t)\tilde{b}(t)-\hat{q}^{0}%
(t)\tilde{\sigma}(t)-\langle\hat{p}^{1}(t),\tilde{M}^{\prime}(t)\rangle\}dt],
\end{array}
\label{i1}%
\end{equation}
where $\tilde{b}(t)=\check{b}(t)-\hat{b}(t)$ etc. By the concavity of $g$ and
the terminal values of the BSDEs \eqref{eqp0}, \eqref{eqp1}, we have%
\[%
\begin{array}
[c]{lll}%
I_{2} & \leq\mathbb{E}[\tfrac{\partial g}{\partial x}(T)\tilde{X}%
(T)+\langle\nabla_{m}g(T),\tilde{M}(T)\rangle] & =\mathbb{E}[\hat{p}%
^{0}(T)\tilde{X}(T)+\langle\hat{p}^{1}(T),\tilde{M}(T)\rangle].
\end{array}
\]
Applying the It\^{o} formula to $\hat{p}^{0}(t)\tilde{X}(t)$ and $\langle
\hat{p}^{1}(t),\tilde{M}(t)\rangle$, we get
\begin{align}
I_{2} &  \leq\mathbb{E}[\hat{p}^{0}(T)\tilde{X}(T)+\langle\hat{p}%
^{1}(T),\tilde{M}(T)\rangle]\nonumber\\
&  =\mathbb{E}[%
{\textstyle\int_{0}^{T}}
\hat{p}^{0}(t)d\tilde{X}(t)+%
{\textstyle\int_{0}^{T}}
\tilde{X}(t)d\hat{p}^{0}(t)+%
{\textstyle\int_{0}^{T}}
\hat{q}^{0}(t)\tilde{\sigma}(t)dt\nonumber\\
&  +\mathbb{E}[%
{\textstyle\int_{0}^{T}}
\langle\hat{p}^{{1}}(t),d\tilde{M}(t)\rangle+%
{\textstyle\int_{0}^{T}}
\tilde{M}(t)d\hat{p}^{1}(t)]\nonumber\\
&  =\mathbb{E}[{\textstyle\int_{0}^{T}}\hat{p}^{0}(t)\tilde{b}%
(t)dt-{\textstyle\int_{0}^{T}}\tfrac{\partial\hat{H}_{0}}{\partial x}%
(t)\tilde{X}(t)dt-{\textstyle\int_{0}^{T}}\mathbb{E}[\nabla_{\bar{x}}^{\ast
}\hat{H}_{0}(t)|\mathcal{F}_{t}]\widetilde{X}(t)dt\nonumber\\
&  -{\int_{0}^{T}}\tfrac{\partial\hat{h}}{\partial x}(t)\widetilde{X}%
(t)d\hat{\xi}(t)+{\textstyle\int_{0}^{T}}\hat{q}^{0}(t)\tilde{\sigma
}(t)dt+{\textstyle\int_{0}^{T}}\langle\hat{p}^{1}(t),\tilde{M}^{\prime
}(t)\rangle dt\nonumber\\
&  -{\textstyle\int_{0}^{T}}\langle\nabla_{m}\hat{H}_{0}(t),\tilde
{M}(t)\rangle dt-{\textstyle\int_{0}^{T}}\mathbb{E}[\nabla_{\bar{m}}^{\ast
}\hat{H}_{0}(t)|\mathcal{F}_{t}]\tilde{M}(t)dt],\label{I2}%
\end{align}
where we have used that the $dB(t)$ and $\tilde{N}(dt,d\zeta)$ integrals with
the necessary integrability property are martingales and then have mean zero.
Substituting $\left(  \ref{i1}\right)  $ and $\left(  \ref{I2}\right)  $ in
\eqref{eq6.16}, yields
\begin{align*}
&  J(\xi)-J(\hat{\xi})\\
&  \leq\mathbb{E}[%
{\textstyle\int_{0}^{T}}
\{H_{0}(t)-\hat{H}_{0}(t)-\tfrac{\partial\hat{H}_{0}}{\partial x}(t)\tilde
{X}(t)-\langle\nabla_{\bar{x}}\hat{H}_{0}(t),\tilde{X}_{t}\rangle\\
&  -\langle\nabla_{m}\hat{H}_{0}(t),\tilde{M}(t)\rangle-\langle\nabla_{\bar
{m}}\hat{H}_{0}(t),\tilde{M}_{t}\rangle\}dt+{\textstyle\int_{0}^{T}}%
h(t)d\xi(t)\\
&  -{\textstyle\int_{0}^{T}}\hat{h}(t)d\hat{\xi}(t)-{\textstyle\int_{0}^{T}%
}\tfrac{\partial\hat{h}}{\partial x}(t)\widetilde{X}(t)d\hat{\xi}(t)\\
&  +{\textstyle\int_{0}^{T}}(\lambda(t)\hat{p}^{0}(t)+h(t))d\xi
(t)-{\textstyle\int_{0}^{T}}(\lambda(t)\hat{p}^{0}(t)+\hat{h}(t))d\hat{\xi
}(t)\\
&  -{\textstyle\int_{0}^{T}}(\lambda(t)\hat{p}^{0}(t)+h(t))d\xi
(t)+{\textstyle\int_{0}^{T}}(\lambda(t)\hat{p}^{0}(t)+\hat{h}(t))d\hat{\xi
}(t)].
\end{align*}
By the concavity of $dH$ and the fact that the process $\xi$ is $\mathbb{G}%
$-adapted, we obtain%
\begin{align}
J(\xi)-J(\hat{\xi}) &  \leq\mathbb{E}[%
{\textstyle\int_{0}^{T}}
\tfrac{\partial\hat{H}_{0}}{\partial\xi}(t)(\xi(t)-\hat{\xi}%
(t))dt+{\textstyle\int_{0}^{T}}(\lambda(t)\hat{p}^{0}(t)+h(t)(d\xi
(t)-d\hat{\xi}(t))]\nonumber\\
&  =\mathbb{E}[%
{\textstyle\int_{0}^{T}}
\mathbb{E}(\tfrac{\partial\hat{H}_{0}}{\partial\xi}(t)(\xi(t)-\hat{\xi
}(t))+\hat{h}(t)(d\xi(t)-d\hat{\xi}(t))|\mathcal{G}_{t})dt]\nonumber\\
&  =\mathbb{E}[{\textstyle\int_{0}^{T}}\langle\mathbb{E}(\nabla_{\xi}\hat
{H}(t)|\mathcal{G}_{t}),\xi(t)-\hat{\xi}(t)\rangle dt]\leq0,\nonumber
\end{align}
where $\frac{\partial\hat{H}_{0}}{\partial\xi}=\nabla_{\xi}\hat{H}_{0}.$ The
last equality holds because $\xi=\hat{\xi}$ maximizes the random measure
$dH(t,\hat{X}(t),\hat{X}_{t},\hat{M}(t),\hat{M}_{t},\xi,\hat{p}^{0}(t),\hat
{q}^{0}(t),\hat{p}^{1}(t))$ at $\xi=\hat{\xi}$. \fproof

\noindent From the above result, we can deduce the following \emph{sufficient
variational inequalities}.

\begin{theorem}
[\emph{Sufficient variational inequalities}]Suppose that $H_{0}$ does not
depend on $\xi$,i.e. that
\[
\tfrac{\partial H_{0}}{\partial\xi}=0,
\]
and that the following variational inequalities hold:
\begin{align}
&  (i)\quad\EE[\lambda(t) \hat{p}^{0}(t)+h(t,\hat{X}(t)) | \mathcal{G}_t]\leq
0,\label{eq5.28}\\
&  (ii)\quad
\EE[\lambda(t) \hat{p}^{0}(t)+h(t,\hat{X}(t)) | \mathcal{G}_t]d\hat{\xi}(t)=0.
\label{eq5.29}%
\end{align}
Then $\hat{\xi}$ is an optimal singular control.\newline
\end{theorem}

\dproof Suppose \eqref{eq5.28} - \eqref{eq5.29} hold. Then for $\xi\in\Xi$ we
have
\[
\EE[\lambda(t) \hat{p}^{0}(t)+h(t,\hat{X}(t)) | \mathcal{G}_t]d\xi
(t)\leq0=\EE[\lambda(t) \hat{p}^{0}(t)+h(t,\hat{X}(t)) | \mathcal{G}_t]d\hat
{\xi}(t).
\]
Since $H_{0}$ does not depend on $\xi$, it follows that \eqref{eq5.21} hold.
\fproof\newline

\subsection{A necessary maximum principle for singular mean-field control}

\label{sec5}

In the previous section we gave a verification theorem, stating that if a
given control $\hat{\xi}$ satisfies \eqref{eq3.10a}-\eqref{eq5.21}, then it is
indeed optimal for the singular mean-field control problem. We now establish a
partial converse, implying that if a control $\hat{\xi}$ is optimal for the
singular mean-field control problem, then it is a conditional critical point
for the Hamiltonian. \newline For $\xi\in\Xi$, let $\mathcal{V}(\xi)$ denote
the set of $\mathbb{G}$-adapted processes $\eta$ of finite variation such that
there exists $\varepsilon=\varepsilon(\xi)>0$ satisfying
\begin{equation}
\xi+a\eta\in\Xi\text{ for all }a\in\lbrack0,\varepsilon].\label{eq5.1a}%
\end{equation}
Note that the following processes $\eta_{i}(s),i=1,2,3$ belong to
$\mathcal{V}(\xi)$:
\begin{align*}
\eta_{1}(s) &  :=\alpha(\omega)\chi_{\lbrack t,T]}(s),\text{ where }%
t\in\lbrack0,T],\alpha>0\text{ is }\mathcal{G}_{t}\text{-measurable },\\
\eta_{2}(s) &  :=\xi(s),\\
\eta_{3}(s) &  :=-\xi(s),s\in\lbrack0,T].
\end{align*}
Then for $\xi\in\Xi$ and $\eta\in\mathcal{V}(\xi)$ we have, by our smoothness
assumptions on the coefficients,
\begin{align}
&  \lim_{a\rightarrow0^{+}}\tfrac{1}{a}(J(\xi+a\eta)-J(\xi))\label{eq5.2}\\
&  ={\mathbb{E}}[%
{\textstyle\int_{0}^{T}}
\{\tfrac{\partial f}{\partial x}(t)Z(t)+\left\langle \nabla_{\bar{x}%
}f(t),Z_{t}\right\rangle +\left\langle \nabla_{m}f(t),DM(t)\right\rangle
\nonumber\\
&  +\left\langle \nabla_{\bar{m}}f(t),DM_{t}\right\rangle \}dt+\tfrac{\partial
f}{\partial\xi}(t)\eta(t)+\tfrac{\partial g}{\partial x}(T)Z(T)\nonumber\\
&  +\left\langle \nabla_{m}g(T),DM(T)\right\rangle +%
{\textstyle\int_{0}^{T}}
\tfrac{\partial h}{\partial x}(t)Z(t)d\xi(t)+%
{\textstyle\int_{0}^{T}}
h(t)d\eta(t)],\nonumber
\end{align}
where
\begin{equation}%
\begin{array}
[c]{l}%
Z(t):=Z_{\eta}(t):=\lim_{a\rightarrow0^{+}}\tfrac{1}{a}(X^{(\xi+a\eta
)}(t)-X^{(\xi)}(t))\\
Z_{t}:=Z_{t,\eta}:=\lim_{a\rightarrow0^{+}}\tfrac{1}{a}(X_{t}^{(\xi+a\eta
)}-X_{t}^{(\xi)})
\end{array}
\label{eq5.3}%
\end{equation}
and
\begin{equation}%
\begin{array}
[c]{l}%
DM(t):=D_{\eta}M(t):=\lim_{a\rightarrow0^{+}}\tfrac{1}{a}(M^{(\xi+a\eta
)}(t)-M^{(\xi)}(t)),\\
DM_{t}:=D_{\eta}M_{t}:=\lim_{a\rightarrow0^{+}}\tfrac{1}{a}(M_{t}^{(\xi
+a\eta)}-M_{t}^{(\xi)}).
\end{array}
\label{eq5.4a}%
\end{equation}
Then
\[
\left\{
\begin{array}
[c]{ll}%
dZ(t) & =[\tfrac{\partial b}{\partial x}(t)Z(t)+\left\langle \nabla_{\bar{x}%
}b(t),Z_{t}\right\rangle +\left\langle \nabla_{m}b(t),DM(t)\right\rangle
+\left\langle \nabla_{\bar{m}}b(t),DM_{t}\right\rangle \\
& +\tfrac{\partial b}{\partial\xi}(t)\eta(t)]dt+[\tfrac{\partial\sigma
}{\partial x}(t)Z(t)+\left\langle \nabla_{\bar{x}}\sigma(t),Z_{t}\right\rangle
+\left\langle \nabla_{m}\sigma(t),DM(t)\right\rangle \\
& +\left\langle \nabla_{\bar{m}}\sigma(t),DM_{t}\right\rangle +\tfrac{\partial
b}{\partial\xi}(t)\eta(t)]dB(t)+\lambda(t)d\eta(t)\;;\\
Z(0) & =0\,,
\end{array}
\right.
\]
and similarly with $dZ_{t},dDM(t)$ and $dDM_{t}$.

\noindent We first state and prove a basic step towards a necessary maximum
principle.

\begin{proposition}
\label{prop}\ Let $\xi\in\Xi$ and choose $\eta\in\mathcal{V}(\xi)$.Then%

\begin{equation}
\tfrac{d}{da}J(\xi+a\eta)|_{a=0}={\mathbb{E}}[%
{\textstyle\int_{0}^{T}}
\tfrac{\partial H_{0}}{\partial\xi}(t)\eta(t)dt+%
{\textstyle\int_{0}^{T}}
\{\lambda(t)p^{0}(t)+h(t)\}d\eta(t)]. \label{eq6.44}%
\end{equation}

\end{proposition}

\noindent{Proof.} \quad Let $\xi\in\Xi$ and $\eta\in\mathcal{V}(\xi)$. Then we
can write
\begin{equation}
\tfrac{d}{da}J(\xi+a\eta)|_{a=0}=A_{1}+A_{2}+A_{3}+A_{4}, \label{eq6.45}%
\end{equation}
where
\begin{align*}
A_{1}  &  ={\mathbb{E}}[%
{\textstyle\int_{0}^{T}}
\{\tfrac{\partial f}{\partial x}(t)Z(t)+\left\langle \nabla_{\bar{x}%
}f(t),Z_{t}\right\rangle +\left\langle \nabla_{m}f(t),DM(t)\right\rangle
+\left\langle \nabla_{\bar{m}}f(t),DM_{t}\right\rangle \}dt],\\
A_{2}  &  ={\mathbb{E}}[%
{\textstyle\int_{0}^{T}}
\tfrac{\partial f}{\partial\xi}(t)\eta(t)dt],\\
A_{3}  &  ={\mathbb{E}}[\tfrac{\partial g}{\partial x}(T)Z(T)+\left\langle
\nabla_{m}g(T),DM(T)\right\rangle ]\\
A_{4}  &  ={\mathbb{E}}[%
{\textstyle\int_{0}^{T}}
\tfrac{\partial h}{\partial x}(t)Z(t)d\xi(t)+h(t)d\eta(t)].
\end{align*}
By the definition of $H_{0}$ we have
\begin{align}
A_{1}  &  ={\mathbb{E}}[%
{\textstyle\int_{0}^{T}}
Z(t)\{\tfrac{\partial H_{0}}{\partial x}(t)-\tfrac{\partial b}{\partial
x}(t)p^{0}(t)-\tfrac{\partial\sigma}{\partial x}(t)q^{0}(t)\}dt \label{eq5.20}%
\\
&  +%
{\textstyle\int_{0}^{T}}
\left\langle \nabla_{\bar{x}}H_{0}(t)-\nabla_{\bar{x}}b(t)p^{0}(t)-\nabla
_{\bar{x}}\sigma(t)q^{0}(t),Z_{t}\right\rangle dt\nonumber\\
&  +%
{\textstyle\int_{0}^{T}}
\left\langle \nabla_{m}H_{0}(t)-\nabla_{m}b(t)p^{0}(t)-\nabla_{m}%
\sigma(t)q^{0}(t),DM(t)\right\rangle dt\nonumber\\
&  +%
{\textstyle\int_{0}^{T}}
\left\langle \nabla_{\bar{m}}H_{0}(t)-\nabla_{\bar{m}}b(t)p^{0}(t)-\nabla
_{\bar{m}}\sigma(t)q^{0}(t),DM_{t}\right\rangle \}dt],\nonumber
\end{align}
and
\[
A_{2}={\mathbb{E}}[%
{\textstyle\int_{0}^{T}}
\{\tfrac{\partial H_{0}}{\partial\xi}(t)-\tfrac{\partial b}{\partial\xi
}(t)p^{0}(t)-\tfrac{\partial\sigma}{\partial\xi}(t)q^{0}(t)\}\eta(t)dt].
\]
By the terminal conditions of $p^{0}(T)$, $p^{1}(T)$ (see
\eqref{eqp0}-\eqref{eqp1}) and the It\^{o} formula, we have
\begin{align}
A_{3}  &  ={\mathbb{E}}[p^{0}(T)Z(T)+\left\langle p^{1}(T),DM(T)\right\rangle
]\label{eq6.49}\\
\text{ }  &  ={\mathbb{E}}[%
{\textstyle\int_{0}^{T}}
p^{0}(t)dZ(t)+%
{\textstyle\int_{0}^{T}}
Z(t)dp^{0}(t)\nonumber\\
&  +%
{\textstyle\int_{0}^{T}}
q^{0}(t)\{\tfrac{\partial\sigma}{\partial x}(t)Z(t)+\left\langle \nabla
_{\bar{x}}\sigma(t),Z(t)\right\rangle +\left\langle \nabla_{m}\sigma
(t),DM(t)\right\rangle \nonumber\\
&  +\left\langle \nabla_{\bar{m}}\sigma(t),DM(t)\right\rangle +\tfrac
{\partial\sigma}{\partial\xi}(t)\eta(t)\}dt\nonumber\\
&  +\left\langle p^{1}(t),dDM(t)\right\rangle +\left\langle DM(t),dp^{1}%
(t)\right\rangle \nonumber\\
&  ={\mathbb{E}}[%
{\textstyle\int_{0}^{T}}
p^{0}(t)\{\tfrac{\partial b}{\partial x}(t)Z(t)+\left\langle \nabla_{\bar{x}%
}b(t),Z_{t}\right\rangle +\left\langle \nabla_{m}b(t),DM(t)\right\rangle
\nonumber\\
&  \text{ }+\left\langle \nabla_{\bar{m}}b(t),DM_{t}\right\rangle
+\tfrac{\partial b}{\partial\xi}(t)\eta(t)\}dt\nonumber\\
&  +%
{\textstyle\int_{0}^{T}}
q^{0}(t)\{\tfrac{\partial\sigma}{\partial x}(t)Z(t)+\left\langle \nabla
_{\bar{x}}\sigma(t),Z_{t}\right\rangle +\left\langle \nabla_{m}\sigma
(t),DM(t)\right\rangle \nonumber\\
&  \text{ }+\left\langle \nabla_{\bar{m}}\sigma(t),DM_{t}\right\rangle
+\tfrac{\partial\sigma}{\partial\xi}(t)\eta(t)\}dt\nonumber\\
&  +{\textstyle\int_{0}^{T}}p^{0}(t)\lambda(t)d\eta(t)+{\textstyle\int_{0}%
^{T}}\big\{Z(t)(-\{\tfrac{\partial H_{0}}{\partial x}(t)+{\mathbb{E}}%
(\nabla_{\bar{x}}^{\ast}H_{0}(t)|\mathcal{F}_{t})\})\nonumber\\
&  -\left\langle \nabla_{m}H_{0}(t)+\mathbb{E}[\nabla_{\bar{m}}^{\ast}%
H_{0}(t)|\mathcal{F}_{t}],DM(t)\right\rangle \big\}dt-{\textstyle\int_{0}^{T}%
}\tfrac{\partial h}{\partial x}(t)Z(t)d\xi(t)].\nonumber
\end{align}
Combining \eqref{eq6.45}-\eqref{eq6.49} and using \eqref{eq6.7a}, we get
\eqref{eq6.44}. \hfill$\square$ \bigskip

\begin{theorem}
[Necessary maximum principle for mean-field singular control]\label{th5.1b}%

Suppose $\hat{\xi}\in\Xi$ is optimal, i.e. satisfies \eqref{eq6.4}. Suppose
that
\[
\tfrac{\partial H_{0}}{\partial\xi}=0.
\]
Then the following variational inequalities hold:
\begin{align}
&  (i)\quad\EE[\lambda(t) \hat{p}^0(t) + h(t)|\mathcal{G}_t ]\leq0\text{ for
all }t\in\lbrack0,T]\text{ a.s. }\quad\text{and}\label{eq1.17a}\\
&  (ii)\quad\EE[\lambda(t)\hat{p}^0(t)+\hat{h}(t)|\mathcal{G}_{t}]d\hat{\xi
}(t)=0\text{ for all }t\in\lbrack0,T]\,\text{ a.s. } \label{eq1.17b}%
\end{align}

\end{theorem}

\noindent{Proof.} \quad From Proposition (\ref{prop}) we have, since $\hat
{\xi}$ is optimal,%

\begin{equation}
0\geq\tfrac{d}{da}J(\hat{\xi}+a\eta)|_{a=0}={\mathbb{E}}[%
{\textstyle\int_{0}^{T}}
\{\lambda(t)\hat{p}^{0}(t)+\hat{h}(t)\}d\eta(t)], \label{eq5.17a}%
\end{equation}
for all $\eta\in\mathcal{V}(\hat{\xi})$.\newline If we choose $\eta$ to be a
pure jump process of the form
\[
\eta(s)=%
{\textstyle\sum_{0<t_{i}\leq s}}
\alpha(t_{i}),
\]
where $\alpha(s)>0$ is $\mathcal{G}_{s}$-measurable for all $s$, then $\eta
\in\mathcal{V}(\hat{\xi})$ and \eqref{eq5.17a} gives
\[
{\mathbb{E}}[\{\lambda(t)\hat{p}^{0}(t)+\hat{h}(t)\}\alpha(t_{i})]\leq0\text{
for all }t_{i}\text{ a.s. }%
\]
Since this holds for all such $\eta$ with arbitrary $t_{i}$, we conclude that
\begin{equation}
{\mathbb{E}}[\lambda(t)\hat{p}^{0}(t)+\hat{h}(t)|\mathcal{G}_{t}]\leq0\text{
for all }t\in\lbrack0,T]\text{ a.s. } \label{eq5.23a}%
\end{equation}
Finally, applying \eqref{eq5.17a} to $\eta_{1}:=\hat{\xi}\in\mathcal{V}%
(\hat{\xi})$ and then to $\eta_{2}:=\hat{\xi}\in\mathcal{V}(\hat{\xi})$ we
get, for all $t\in\lbrack0,T]$,
\begin{equation}
{\mathbb{E}}[\lambda(t)\hat{p}^{0}(t)+\hat{h}(t)|\mathcal{G}_{t}]d\hat{\xi
}(t)=0\text{ for all }t\in\lbrack0,T]\text{ a.s. } \label{eq5.24a}%
\end{equation}
With \eqref{eq5.23a} and \eqref{eq5.24a} the proof is complete. \hfill
$\square$ \bigskip

\section{Application to optimal stopping}

\label{sec6.1} From now on, let us assume, in addition to
$$\frac{\partial H_0}{\partial \xi}=0,$$
 that
\begin{align}
 \lambda(t)&=-\lambda_{0} \text{ where } \lambda_{0} > 0, \text { and } \\
 \mathbb{G}&=\mathbb{F}.
\end{align}

\noindent Then, dividing by $\lambda_0$ in \eqref{eq1.17a} - \eqref{eq1.17b} we get
\begin{align}
\label{eq6.2} &  (i) \quad \hat{p}^0(t) \geq \frac{1}{\lambda_{0}} \hat {h}(t)) \text{ for all } t\in[0, T] \text{ a.s. } \quad\text{and}\\
&  (ii) \quad\Big\{\hat{p}^0(t) - \frac{1}{\lambda_{0}} \hat{h}(t)\Big\} d\hat{\xi}(t) = 0 \text{ for all } t\in[0, T]\, \text{ a.s. } \label{eq6.3}%
\end{align}

\noindent Comparing with \eqref{eq3.1}, we see that \eqref{eq6.2}-\eqref{eq6.3},
together with the singular BSDE \eqref{eqp0} for ${p}^0=\hat{p}^0,{q}^0=\hat{q}^0,\xi=\hat{\xi},$
constitutes an AMBSDEs related to the type discussed
in Section 3 above, with \newline%
\begin{equation} \label{eq5.48a}
S(t)=\frac{1}{\lambda_{0}} \hat{h}(t),
\end{equation}
and
\begin{align}
Y(t)&:= \hat{p}^0(t)\label{eq5.48} ,\\
Z(t)&:=\hat{q}^0(t),\\
dK(t)&:= \frac{\partial\hat{h}}{\partial x}(t) d\hat{\xi}(t). \label{eq5.49}
\end{align}

\noindent We summarize what we have proved as follows:

\begin{theorem}
Suppose $\hat{\xi}$ is an optimal control for the singular control problem
\eqref{eq6.1} - \eqref{eq6.4}, with corresponding optimal processes $\hat
{X}(t),\hat{X}_{t},\hat{M}(t),\hat{M}_{t}$. Define $S, Y,Z,K$ as in \eqref{eq5.48a}, \eqref{eq5.48}, \eqref{eq5.49}.
\
Then $\hat{X}$ together with $(Y,Z,K)$
solve the following \emph{forward-backward memory-advanced mean-field singular reflected
system}:
\end{theorem}

\begin{itemize}
\item (i) Forward mean-field \emph{memory} singular SDE in $\hat{X}$:

\begin{equation}
\left\{
\begin{array}
[c]{l}%
d\hat{X}(t)=b(t,\hat{X}(t),\hat{X}_{t},\hat{M}(t),\hat{M}_{t})dt \\
+\sigma(t,\hat{X}(t),\hat{X}_{t},\hat{M}(t),\hat{M}_{t})dB(t)
-\lambda_{0}d\hat{\xi
}(t);\quad t\in\lbrack0,T]\\
X(t)=\alpha(t),\quad t\in\lbrack-\delta,0],
\end{array}
\right.  %
\end{equation}

\item (ii) \emph{Advanced} reflected BSDE in $(Y,Z,K)$ (for given $\hat
{X}(t)$):
\[%
\begin{cases}
& dY(t)=-\big\{\frac{\partial\hat{H}_{0}}{\partial x}(t)+{\mathbb{E}%
}[\nabla^{*}_{\bar{x}}\hat{H}_{0}(t)|\mathcal{F}_{t}]\big \}dt\\
& -dK(t)+Z(t)dB(t);\quad t\in\lbrack0,T],\\ & Y(t)\geq S(t);\quad t\in\lbrack0,T],\\
& [Y(t)-S(t)]dK(t)=0;\quad t\in\lbrack0,T],\\
& Y(T)=\frac{\partial g}{\partial x}(T).
\end{cases}
\]

\end{itemize}

\subsection{Connection to optimal stopping of memory mean-field SDE}

\vskip 0.3cm
If we combine the results above, we get
\begin{theorem}
Suppose $\hat{\xi}$ is an optimal control for the singular control problem
\eqref{eq6.1} - \eqref{eq6.4}, with corresponding optimal processes $\hat
{X}(t),\hat{X}_{t},\hat{M}(t),\hat{M}_{t}$ and adjoint processes $\hat{p}^{0}(t),\hat{q}^{0}(t)$.
Put
\begin{equation}
R=\frac{\partial g}{\partial x}(T).
\end{equation}
Let
$$S(t), (Y(t),Z(t),K(t))$$
 be as above and define
 \begin{align}
 F(t)&:=F(t,\hat{X}(t),\hat{M}(t),\hat{X}_t,\hat{M}_t, Y(t),Z(t),Y^{t},\mathcal{Z}^{t})\nonumber\\
 &:=\frac{\partial\hat{H}_{0}}{\partial x}(t)+{\mathbb{E}%
}[\nabla^{*}_{\bar{x}}\hat{H}_{0}(t)|\mathcal{F}_{t}].
 \end{align}
\begin{description}
\item[(i)] Then, for each $t\in\left[  0,T\right], Y(t)$ is the solution of the optimal stopping problem%
\begin{align}
Y(t)=\underset{\tau\in\mathcal{T}_{[t,T]}}{ess\sup}\Big\{ \mathbb{E} [%
{\textstyle\int_{t}^{\tau}}
F(s)ds+S(\tau)\mathbf{1}_{\tau<T}%
+R\mathbf{1}_{\tau=T}|\mathcal{F}_{t}]\Big\}.
\end{align}

\item[(ii)] Moreover, for $t \in [0,T]$  the solution process $K(t)$ is given by%
\begin{align}
&K(T)-K(T-t)\nonumber\\
&=\underset{s\leq t}{\max}\Big\{R+%
\int_{T-s}^{T}
F(r)dr-\int _{T-s}^{T} Z(r)dB(r)-S(T-s)\Big\}^{-},
\end{align}
where $x^{-}=\max(-x,0),$ and an optimal stopping time $\hat{\tau}_{t}$ is given by%
\begin{align*}
\hat{\tau}_{t}:&=\inf\{s\in\lbrack t,T],Y(s)\leq S(s)\}\wedge T\\
&=\inf\{s\in\lbrack t,T],K(s) > K(t)\}\wedge T.
\end{align*}

\item[(iii)] In particular, if we choose $t=0$ we get that
\begin{align*}
\hat{\tau}_{0}:&=\inf\{s\in\lbrack 0,T],Y(s)\leq S(s)\}\wedge T \\
&=\inf\{s\in\lbrack 0,T],K(s) > 0\}\wedge T
\end{align*}
solves the optimal stopping problem
\[
Y(0)=\sup_{\tau\in\mathcal{T}_{[0,T]}}\mathbb{E}[{\textstyle\int_{0}^{\tau}}
F(s)ds+S(\tau)\mathbf{1}_{\tau<T}%
+R\mathbf{1}_{\tau=T}] .
\]

\end{description}
\end{theorem}

\end{document}